\documentclass[a4paper,10pt]{article}
\usepackage[utf8]{inputenc}

\usepackage{amsmath}
\usepackage{amssymb}
\usepackage{siunitx}

\usepackage{graphicx}

\usepackage{geometry}
\newcommand{\width}[1]{{\text{\rm w}{\left(#1\right)}}}
\newcommand{\Mid}{\text{\rm m}}
\newcommand{\Inf}[1]{{\underline{#1}}}
\newcommand{\Sup}[1]{{\overline{#1}}}
\newcommand{\Rset}{\mathbf{\mathbb{R}}}
\newcommand{\IRset}{\mathbf{\mathbb{IR}}}

\usepackage{amsthm}

\newtheorem{proposition}{Proposition}

\newtheorem{definition}{Definition}

\newtheorem{example}{Example}

\begin{document}

\title{Sensitivity-based Heuristic for Guaranteed Global Optimization with Nonlinear Ordinary Differential Equations
}


\author{Julien Alexandre dit Sandretto
}

%

\date{May 2018}

\maketitle

\begin{abstract}
We focus on interval algorithms for computing guaranteed enclosures of the solutions of constrained global optimization problems where differential constraints occur. 
To solve such a problem of global optimization with nonlinear ordinary differential equations, a branch and bound algorithm can be used based on guaranteed numerical integration methods. 
Nevertheless, this kind of algorithms is expensive in term of computation. Defining new methods to reduce the number of  branches is still a challenge. 
Bisection based on the smear value is known to be often the most efficient heuristic for branching algorithms. 
This heuristic consists in bisecting in the coordinate direction for which the values of the considered function change the most ``rapidly''.

We propose to define a smear-like function using the sensitivity function obtained from the differentiation of ordinary differential equation with respect to parameters. 
The sensitivity has been already used in validated simulation for local optimization but not as a bisection heuristic. 
We implement this heuristic in a branch and bound algorithm to solve a problem of global optimization with nonlinear ordinary differential equations. 
Experiments show that the gain in term of number of branches could be up to 30\%.

\end{abstract}

\section{Introduction}
A physical system, such as a physical phenomenon, a chemical reaction, or more recently an autonomous vehicle, is often modeled by an Ordinary Differential Equation (ODE), or eventually a Differential Algebraic Equation (DAE). These models are parametrized to represent as accurately as possible the actual system. Optimization of such models are an important but difficult problem. Obtaining the parameters such that a system is optimal can be addressed in optimal design, calibration or optimal control \cite{rauh2009}. 
Three main approaches are classically used to solve such a problem: the calculus of variations in their infinite dimensional spaces \cite{pontryagin_1987}, by discretization of time, state variables and parameters (Monte Carlo), or complete collocation \cite{neuman_1973,esposito_2000}. In general, these approaches provide local minima (except for certain convex problems), and finding a global minimum is a challenge. 

For a rigorous global optimization approach, the branch-and-bound algorithm has shown its efficiency \cite{papamichail2002rigorous}. It consists in a deterministic spatial branch-and-bound to solve constrained optimization systems, considering the differential equation as a constraint. Interval arithmetic is used in \cite{papamichail2002rigorous} like in this paper. In \cite{Singer2006}, authors present a method for global optimization with nonlinear ODEs based on a branch-and-bound algorithm exploiting some constructed relaxations (convex and concave). The introduction of \cite{Singer2006} proposes a detailed state of the art. 
Interval analysis has been used in other approaches than the one explained in \cite{papamichail2002rigorous}, such as in \cite{lin_2006} and \cite{rauh2013interval}. This latter presents a method for a local search based on sensitivity analysis. Sensitivity analysis has been used in control of uncertain systems \cite{Rauh2012}. 
We present in this paper a heuristic for bisection also based on sensitivity analysis. It is important to notice that both could be used in a single tool. 

In term of heuristics of bisection, nothing has been clearly done for the particular case of global optimization with ODEs. However, it is common to define the most efficient bisection strategy with respect to a given satisfaction problem solved by the help of a branch-and-prune (or branch-and-contract) \cite{10.1007/978-3-642-33558-7_23}. Indeed, the strategy chosen has an important effect on the complexity of the solver and on the computation time (not on the results obtained). 
The most used heuristics are {\it Round Robin} (bisecting each component in turn) and {\it Largest First} (bisecting the largest component). Another one, the {\it Smear} heuristic consists to bisect in the dimension which has the most important effect. This mode is based on the Smear function as defined by \cite{Kearfott:1990:AIP:78928.78931}, and has been declined in many versions \cite{araya_2013}. 

In this paper, a Smear-like heuristic dedicated to the global optimization with ODEs is proposed to reduce the number of branching in a branch-and-bound algorithm. First, the considered problem is described in Section~\ref{sec:dyn}. In the Section~\ref{sec:preliminaries}, the preliminaries required in this paper are detailed. Then, the branch-and-bound algorithm dedicated to the guaranteed global optimization with ODEs is given in Section~\ref{sec:ggo}. The main contribution is presented in Section~\ref{sec:smear}. Some examples are given in Section~\ref{sec:case} to show the impact of the Smear heuristic in the particular case of dynamic optimization. Finally, we conclude in Section~\ref{sec:concl}.

\section{Dynamic Optimization}
\label{sec:dyn}
This section details the abstract mathematical formulation of the problem considered in this paper. In the rest of the paper, states are denoted by $y$, parameters are notified by $p$, while $t$ is reserved for time ($t\in \Rset^+$). Vector-valued quantities are indicated in bold font.  

\subsection{Problem statement}
\label{problem}

The general problem of global dynamic optimization is stated as follows:

 \begin{equation*}
  \min_{\mathbf{p}} J(\mathbf{p}) = \phi(\mathbf{y}(t_f;\mathbf{p}), \mathbf{p}) + \int_{t_0}^{t_f} g(t, \mathbf{y}(t;\mathbf{p}), \mathbf{p}) dt
 \end{equation*}

  with $\mathbf{y}(t;\mathbf{p})$ solution of the Initial Value Problem (IVP):\\
  \begin{equation*}
  \label{eq:go}
  \left\{
    \begin{aligned}
      \dot{\mathbf{y}}(t;\mathbf{p}) & = \mathbf{f}(t,\mathbf{y}(t;\mathbf{p});\mathbf{p})
     \\
      \mathbf{y}(0) &= \mathbf{y}_0, 
      \quad t \in[t_0, t_{f}]\enspace
  \end{aligned}
\right.
\end{equation*}

$\phi: \mathcal{Y} \times \mathcal{P} \mapsto \Rset$ is a terminal cost, $g: [t_0,t_{f}] \times \mathcal{Y} \times \mathcal{P} \mapsto \Rset$ a continuous cost (Lebesgue integrable) and $\mathbf{f}: \Rset \times \mathcal{Y} \times \mathcal{P} \mapsto \Rset^n$ is continuous in $t$ and Lipschitz and $C^{k}$ ($k$ greater than the order of the integration method) in $\mathbf{y}$ to accept a unique solution $\mathbf{y}(t)$. The sets are defined such that $\mathcal{Y} \subset \Rset^n$ and $\mathcal{P} \subset \Rset^m$.

\subsection{Sensitivity analysis}
\label{sec:adjoint}
The numerical sensitivity analysis of parametrized IVPs is often used in problem involving differential equations. For examples, it can be used to study of dynamical systems, parameter optimization or dynamic optimization. 
The forward sensitivity analysis is the (first order) sensitivity with respect to the model parameter $p_i$ defined as the vector
\begin{equation}
 s_i(t;\mathbf{p}) = \{\frac{\partial }{\partial p_i}\} \mathbf{y}(t;\mathbf{p}), \quad i=1,\dots,m
\end{equation}

Using internal differentiation, the chain rule and the Clairaut's Theorem, the forward sensitivity analysis can be obtained by the numerical integration of the IVP described by
\begin{equation}
 \dot{s}_i = \frac{\partial f}{\partial y} s_i + \frac{\partial f}{\partial p_i}, \quad s_i(t_0) = \frac{\partial y_0(\mathbf{p})}{\partial p_i}, \quad i=1,\dots,m
\end{equation}

This ODE is called the sensitivity equation, while its solution is the forward sensitivity. Sensitivity function is a common tool in the field of dynamical systems for their analysis, parameter identification or optimization, such as for physical systems \cite{hall1983physical,hamby1994review} and more largely for any complex systems \cite{marchuk2013adjoint}. 

\section{Preliminaries}
\label{sec:preliminaries}
In this section, the tools required by the contribution presented in this paper are described.

\subsection{Interval Analysis}
   \label{sec:interval-analysis}

The simplest and most common way to represent and manipulate sets of
values is \textit{interval arithmetic} (see~\cite{Moore66}). An
interval $[x_i]=[\underline{x_i},\overline{x_i}]$ defines the set of
real numbers $x_i$ such that $\underline{x_i} \leq x_i \leq
\overline{x_i}$. $\IRset$ denotes the set of all intervals over
real numbers. The size (or width) of $[x_i]$ is denoted by
$w([x_i])=\overline{x_i}-\underline{x_i}$.

\textit{Interval arithmetic} extends to $\IRset$ elementary functions
over $\Rset$. For instance, the interval sum, \textit{i.e.},
$[x_1]+[x_2]=[\underline{x_1}+\underline{x_2},
\overline{x_1}+\overline{x_2}]$, encloses the image of the sum
function over its arguments.

An interval vector or a \emph{box}
$[\mathbf{x}] \in \IRset^n$, is a Cartesian product of $n$
intervals. The enclosing property basically defines what is called an
{\it interval extension} or an \emph{inclusion function}.

\begin{definition}[Inclusion function]
  Consider a function $f: \Rset^n \rightarrow \Rset^m$, then $[f]\!:\!
  \IRset^n \rightarrow \IRset^m$ is said to be an inclusion function of $f$ to
  intervals if
  \begin{displaymath}
    \forall [\mathbf{x}] \in \IRset^n, \quad [f]([\mathbf{x}])
    \supseteq \{f(\mathbf{x}), \mathbf{x} \in [\mathbf{x}]\}\enspace.
  \end{displaymath}
\end{definition}

It is possible to define inclusion functions for all elementary
functions such as $\times$, $\div$, $\sin$, $\cos$, $\exp$, and so
on. The \emph{natural} inclusion function is the simplest to obtain:
all occurrences of the real variables are replaced by their interval
counterpart and all arithmetic operations are evaluated using interval
arithmetic. More sophisticated inclusion functions such as the
centered form, or the Taylor inclusion function may also be used
(see~\cite{JKDW01} for more details).

\begin{example}[Interval arithmetic]
  A few examples of arithmetic operations between interval values are
  given
  \begin{eqnarray*}
    \lbrack-2,5\rbrack+\lbrack-8,12\rbrack &=& \lbrack-10,17\rbrack \\
    \lbrack-10,17\rbrack-\lbrack-8,12\rbrack=\lbrack-10,17\rbrack+\lbrack-12,8\rbrack &=& \lbrack-22,25\rbrack \\
    \lbrack-10,17\rbrack-\lbrack-2,5\rbrack &=& \lbrack-15,19\rbrack \\
    \frac{\lbrack-2,5\rbrack} {\lbrack-8,12\rbrack} &=& \lbrack-\infty,\infty\rbrack \label{infi}\\
    \frac{\lbrack3,5\rbrack}{\left\lbrack8,12\right\rbrack} &=& \left\lbrack\frac{3}{12},\frac{5}{8}\right\rbrack\\
    \left\lbrack\frac{3}{12},\frac{5}{8}\right\rbrack \times \left\lbrack8,12\right\rbrack &=& \left\lbrack2,\frac{15}{2}\right\rbrack
  \end{eqnarray*}
  In the first example of division, the result is the interval
  containing all the real numbers because denominator contains~$0$.

  As an example of inclusion function, we consider a function $p$
  defined by
  \begin{displaymath}
    p(x,y) = x y + x\enspace.
  \end{displaymath}
  The associated natural inclusion function is
  \begin{displaymath}
    [p]([x],[y]) = [x][y] + [x],
  \end{displaymath}
  in which variables, constants and arithmetic operations have been
  replaced by their interval counterparts. And so \\$p([0,1], [0,1]) =
  [0,2] \subseteq \{ p(x,y) \mid x, y \in [0,1] \} =
  [0,2]$.\hfill$\blacksquare$
\end{example}

\subsection{Validated Numerical Integration Methods}
\label{sec:validated-numerical-integration-methods}

A differential equation is a mathematical relation used to define an
unknown function by its derivatives. These derivatives usually
represent the temporal evolution of a physical quantity. Because the
majority of systems is defined by the evolution of a given quantity,
differential equations are used in engineering, physics, chemistry,
biology or economics as examples. Mathematically, differential
equations have no explicit solutions, except for few particular
cases. Nevertheless, the solution can be numerically approximated with
the help of integration schemes. It is also possible to compute the error introduced 
by these schemes in order to validate them as explained below. Two families of integration schemes have been validated: the Taylor series~\cite{Nedialkov}
and the Runge-Kutta methods \cite{alexandreditsandretto:hal-01243044,alexandreditsandretto:hal-01243053}.

In the following, we consider a \textit{generic} parametric
differential equation as an \emph{interval initial value problem}
(IIVP) defined by
\begin{equation}
  \label{eq:iivp}
  \left\{
    \begin{aligned}
      \dot{\mathbf{y}} & =
      F(t,\mathbf{y},\mathbf{x},\mathbf{p},\mathbf{u})
      \\
      0 & = G(t,\mathbf{y},\mathbf{x},\mathbf{p},\mathbf{u})
      \\
      \mathbf{y}(0) & \in \mathcal{Y}_0, \mathbf{x}(0) \in
      \mathcal{X}_0, \mathbf{p} \in \mathcal{P}, \mathbf{u} \in
      \mathcal{U}, t \in[0, t_{\text{end}}]\enspace,
  \end{aligned}
\right.
\end{equation}
with $F: \Rset \times \Rset^n \times \Rset^m \times \Rset^r \times
\Rset^s \mapsto \Rset^n$ and $G: \Rset \times \Rset^n \times \Rset^m
\times \Rset^r \times \Rset^s \mapsto \Rset^m$. The variable
$\mathbf{y}$ of dimension $n$ is the differential variable while the
variable $\mathbf{x}$ is an algebraic variable of dimension $m$ with
an initial condition $\mathbf{y}(0) \in \mathcal{Y}_0 \subseteq
\Rset^n$ and $\mathbf{x}(0) \in \mathcal{X}_0 \subseteq \Rset^n$. In
other words, differential-algebraic equations (DAE) of index~$1$ are
considered, and in the case of $m=0$, this differential equation
simplifies to an ordinary differential equation (ODE). Note that
usually, the initial values of algebraic variable $\mathbf{x}$ are
computed by numerical algorithms used to solve DAE but we consider it
fixed here for simplicity. Variable~$\mathbf{p} \in \mathcal{P}
\subseteq \Rset^r$ stands for parameters of dimension $r$ and
variable~$\mathbf{u} \in \mathcal{U} \subseteq \Rset^s$ stands for a
control vector of dimension $s$. We assume standard hypotheses on $F$
and $G$ to guarantee the existence and uniqueness of the solution of
such problem.

A validated simulation of a differential equation consists in a
discretization of time, such that $t_0 \leqslant \cdots \leqslant
t_{\text{end}}$, and a computation of enclosures of the set of states
of the system $\mathbf{y}_0$, \dots, $\mathbf{y}_{\text{end}}$, by the
help of a guaranteed integration scheme. In details, a guaranteed
integration scheme is made of
\begin{itemize}
\item an integration method $\Phi(F,G, \mathbf{y}_j,t_j,h)$, starting
  from an initial value $\mathbf{y}_j$ at time $t_j$ and a finite time
  horizon $h$ (the step-size), producing an approximation
  $\mathbf{y}_{j+1}$ at time $t_{j+1} = t_j + h$, of the exact
  solution $\mathbf{y}(t_{j+1}; \mathbf{y}_j)$, that is to say,
  $\mathbf{y}(t_{j+1}; \mathbf{y}_j) \approx \Phi(F,G,
  \mathbf{y}_j,t_j,h)$;
\item a truncation error function
  $\text{lte}_{\Phi}(F,G,\mathbf{y}_j,t_j,h)$, such that
  $\mathbf{y}(t_{j+1}; \mathbf{y}_j) = \Phi(F,G, \mathbf{y}_j,t_j,h) +
  \text{lte}_{\Phi}(F,G,\mathbf{y}_j,t_j,h)$.
\end{itemize}

Basically, a validated numerical integration method is based on a
numerical integration scheme such as Taylor series \cite{Nedialkov} or
Runge-Kutta methods
\cite{alexandreditsandretto:hal-01243044,alexandreditsandretto:hal-01243053}
which is extended with interval analysis tools to bound the
\textit{local truncation error}, the distance between the exact and
the numerical solutions. Mainly, such methods work in two stages at
each integration step, starting from an enclosure $[\mathbf{y}_j] \ni
\mathbf{y}(t_j; \mathbf{y}_0)$ at time $t_j$ of the exact solution, we
proceed by: \renewcommand {\theenumi}{\roman{enumi}}
\begin{enumerate}
\item a computation of an \textit{a priori} enclosure
  $[\widetilde{\mathbf{y}}_{j+1}]$ of the solution $\mathbf{y}(t;
  \mathbf{y}_0)$ for all $t$ in the time interval $[t_j,
  t_{j+1}]$. This stage allows one to prove the existence and the
  uniqueness of the solution.
\item a computation of a tightening of state variable
  $[\mathbf{y}_{j+1}] \ni \mathbf{y}(t_{j+1}; \mathbf{y}_0)$ at time
  $t_{j+1}$ using $[\widetilde{\mathbf{y}}_{j+1}]$ to bound the local truncation error term
  $\text{lte}_{\Phi}(F,G,\mathbf{y}_j,t_j,h)$.
\end{enumerate}

Sometimes, additive constraints can be considered, coming from a
mechanical constraint, an energy conservation, a control law, or
whatever \cite{alexandreditsandretto:hal-01243044}. It is important to
understand that these constraints are linked to the system, valid all
the time, and thus semantically different from a constraint coming
from a satisfaction problem which are valid only at particular time
instants. Finally a \textit{constrained parametric IIVP} is considered
\begin{equation}
  \label{eq:ciivp}
  \left\{
    \begin{aligned}
      \dot{\mathbf{y}} & =
      F(t,\mathbf{y},\mathbf{x},\mathbf{p},\mathbf{u})
      \\
      0 & = G(t,\mathbf{y},\mathbf{x},\mathbf{p},\mathbf{u})
      \\
      0 & = H(t,\mathbf{y},\mathbf{p},\mathbf{u})
      \\
      \mathbf{y}(0) & \in \mathcal{Y}_0, \mathbf{x}(0) \in
      \mathcal{X}_0, \mathbf{p} \in \mathcal{P}, \mathbf{u} \in
      \mathcal{U}, t \in[0, t_{\text{end}}]\enspace,
  \end{aligned}
\right.
\end{equation}
where $F$ and $G$ are defined as in~\eqref{eq:iivp} and $H: \Rset \times \Rset^n \times \Rset^r \times \Rset^s \mapsto
\Rset^c$.

The methods defined in
\cite{alexandreditsandretto:hal-01243044,alexandreditsandretto:hal-01243053}
can take into account such kind of constraints by adapting the first
stage of validated numerical integration methods, that is the
computation of the \textit{a priori} enclosure.

A validated simulation starts with the interval enclosures
$[\mathbf{y}(0)]$, $[\mathbf{x}(0)]$, $[\mathbf{p}]$ and
$[\mathbf{u}]$ of respectively, $\mathcal{Y}_0$, $\mathcal{X}_0$,
$\mathcal{P}$, and $\mathcal{U}$. It produces two lists of boxes:
\begin{itemize}
\item the list of discretization time steps.
  $\{[\mathbf{y}_0],\hdots,[\mathbf{y}_{\text{end}}]\}$;
\item the list of {\it a priori} enclosures:
  $\{[\widetilde{\mathbf{y}}_{0}],\hdots
  ,[\widetilde{\mathbf{y}}_{\text{end}}]\}$.
\end{itemize}
Based on these lists, two functions depending on time can be defined
\begin{equation}
  \label{eq:R}
  R: \left\{\begin{array}{c}
      \Rset \mapsto \IRset^n \\
      t \rightarrow [\mathbf{y}]
    \end{array}\right.
\end{equation}
with $\{ \mathbf{\mathbf{y}}(t; \mathbf{y}_0) : \forall \mathbf{y}_0
\subseteq [\mathbf{y}_0] \} \subseteq [\mathbf{y}]$, and
\begin{equation}
  \label{eq:Rt}
  \widetilde{R}: \left\{\begin{array}{c}
      \IRset \mapsto \IRset^n \\
      \left[\underline{t},\overline{t}\right] \rightarrow [\widetilde{\mathbf{y}}]
    \end{array}\right.
\end{equation}
with $\{ \mathbf{y}(t; \mathbf{y}_0): \forall \mathbf{y}_0 \in
[\mathbf{y}_0] \wedge \forall t \in [\underline{t},\overline{t}]\}
\subseteq [\widetilde{\mathbf{y}}]$.

Function~$R$, defined in \eqref{eq:R}, is obtained by new applications
of validated integration method starting from $[\mathbf{y}_k]$ at
$t_k$, and finishing at $t$ with $t_k < t < t_{k+1}$. Function
$\widetilde{R}$, defined in \eqref{eq:Rt}, is obtained with the union
of $[\widetilde{\mathbf{y}}_k]$ with $k=a,\dots,b$ and $t_a <
\underline{t} < \overline{t} < t_b$. These functions are then strictly
conservative.

More abstractly, the functions $R$ and $\widetilde{R}$ define two
interval enclosures of the solution function of differential equations
defined in Equation~\eqref{eq:ciivp}.

\begin{example}
  We consider an Initial Value Problem based on the following
  constrained differential algebraic equation
  \begin{equation}
    \label{ivp_dae_ctc}
    \left\{\begin{aligned}
        F: & \left\{
          \begin{aligned}
            \dot{y_1} & = -y_3 y_2 - (1+y_3) x_1\\
            \dot{y_2}&= y_3 y_1 -(1+y_3) x_2\\
            \dot{y_3}&=1
          \end{aligned}\right.
          \\[1ex]
        G: & \left\{
        \begin{aligned}
          0 &= (y_1-x_2)/5 - \cos(y_3^2 / 2)\\
          0 &= (y_2-x_1)/5 - \sin(y_3^2 / 2)
        \end{aligned}\right.
        \\[1ex]
        H:& \left\{
        \begin{aligned}
          0 &= x_1^2 + x_2^2 -1
        \end{aligned}
      \right.
    \end{aligned} \right.
  %
  %
\end{equation}
  and the initial values $\mathbf{y}(0)=(5,1,0)$ and
  $\mathbf{x}(0)=(-1,0)$. Interesting variables are $y_1$ and $y_2$,
  the states of the system. A picture of the result of validated
  simulation is given in Figure~\ref{fig:result_val_sim} and
  Figure~\ref{fig:result_val_sim2}.

  \begin{figure}[t]
    \centering
    \includegraphics[width=8cm]{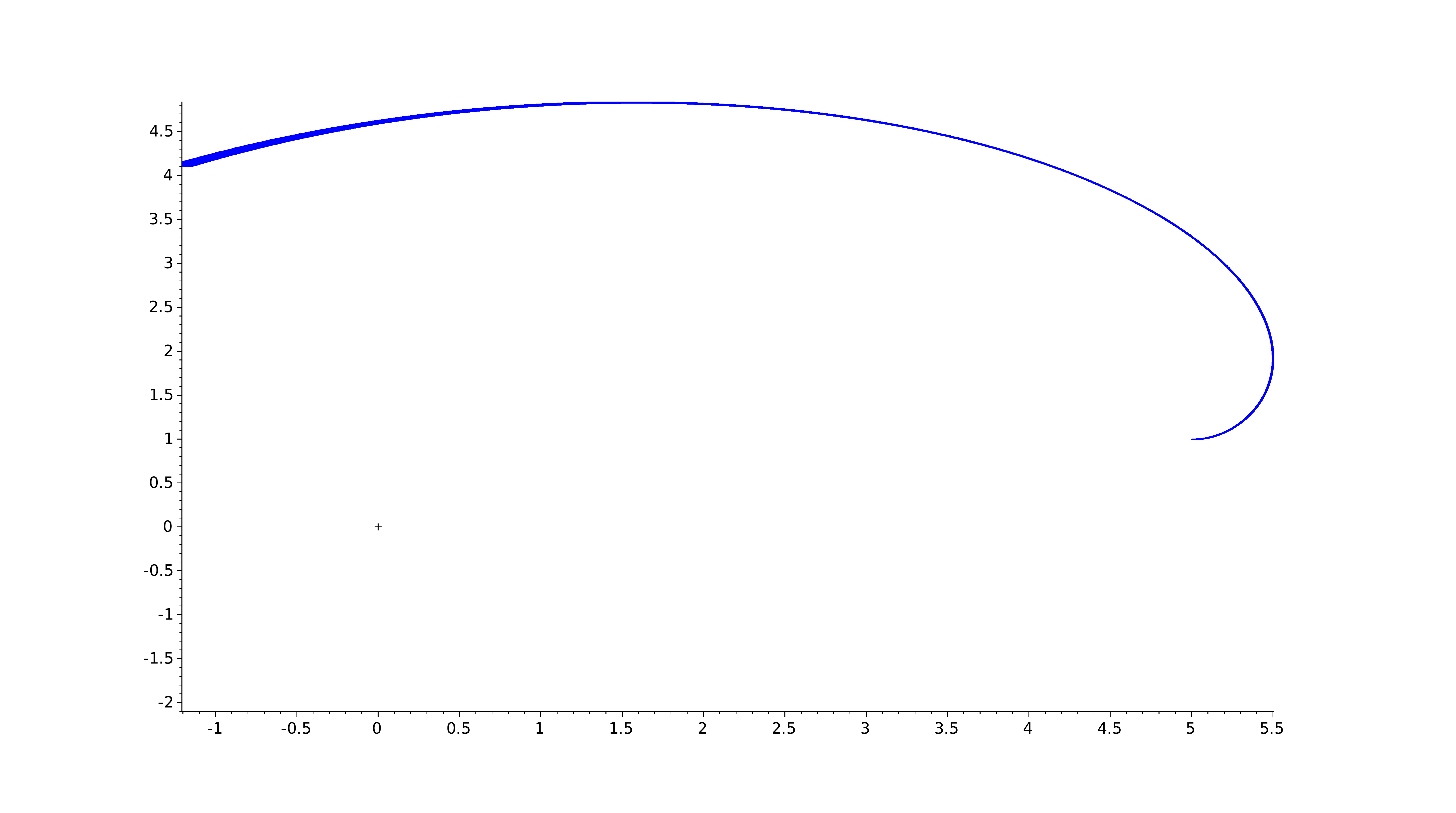}
    \caption{Results of validated simulation: $y_2$ w.r.t. $y_1$.}
    \label{fig:result_val_sim}
  \end{figure}

  \begin{figure}[t]
    \centering
    \includegraphics[width=8cm]{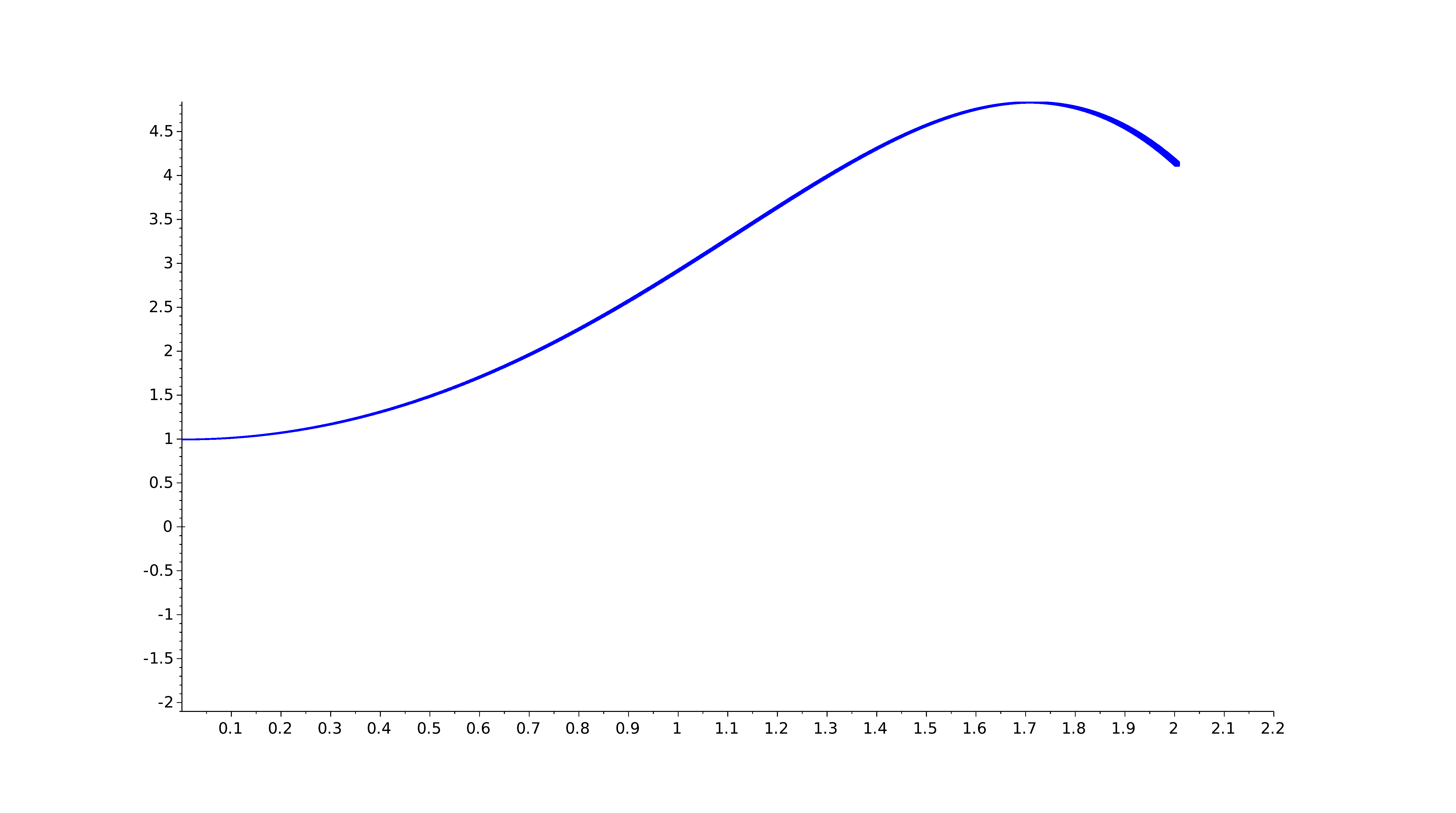}
    \caption{Results of validated simulation: $y_2$ w.r.t. time.}
    \label{fig:result_val_sim2}
  \end{figure}

  It follows an example of the use of the functions providing
  enclosures at given time $R(t)$ and $\widetilde{R}([t])$:
  \begin{itemize}
  \item $R(1.1) = ([5.001, 5.006] ; [3.294, 3.299] ; [1.099, 1.100])$
  \item $\widetilde{R}([0.7,0.8]) =$\\
    $([5.463, 5.495]$; $[1.975, 2.270]$; $[0.699, 0.800])$
  \end{itemize}

  The problem given by Equation~\eqref{ivp_dae_ctc} has a known exact
  solution which is
  \begin{equation*}
    \mathbf{y}(t)=\left\{
      \begin{aligned}
        y_1 & = \sin(t)+5\cos(t^2/2)\\
        y_2 & = \cos(t)+5\sin(t^2/2)\\
        y_3 & = t
      \end{aligned}\right.
    \enspace.
  \end{equation*}

  If we compute the previous enclosures with the exact solution (by
  combining interval arithmetic with bisection method to obtain sharp
  results), we find:
  \begin{itemize}
  \item $\mathbf{y}(1.1) = ([5.003, 5.003] ; [3.297, 3.297] ; [1.1])$
  \item $\mathbf{y}([0.7,0.8]) = ([5.464, 5.494] ; [1.977, 2.269] ; [0.7,0.8])$
  \end{itemize}
  The enclosure property of $R(t)$ and $\widetilde{R}(t)$ is then
  verified.\hfill$\blacksquare$
\end{example}

\section{Guaranteed Global Optimization}
\label{sec:ggo}
The main objective of this paper is to propose a heuristic to improve the interval algorithms for computing guaranteed enclosures of the sets of points which:
\begin{itemize}
 \item validate a set of constraints;
 \item have the global optimal (generally minimal) evaluation of a given cost function.
\end{itemize}

With interval approach, the main idea is to, starting with a box $[\mathbf{x}]$ (supposed to contain the minimum), remove all the infeasible points and the feasible but with a cost greater than the current minimum \cite{hansen_global_2003}. At the fixed point, the resulting box (or list of boxes) contains the constrained global optima. 
If we consider the IVP as a constraint in Problem~\ref{problem}, the interval algorithms for guaranteed global optimization can be applied to solve our problem. Nevertheless, some points need to be clarified before. 

\subsection{Differential Inclusion}
The interval based methods such as the Interval Branch \& Bound algorithm needs to handle intervals in the formulation of the problem. In our case, the branching algorithm is done on the parameter set which has to be considered in its interval form. 

First, we consider the IVP included in problem~\ref{problem}, using interval notation:

\begin{equation}
\label{eq:iivp2}
\begin{aligned}
      \dot{\mathbf{y}}(t;\mathbf{p}) & = \mathbf{f}(t,\mathbf{y}(t;\mathbf{p});\mathbf{p})
     \\
      \text{with }\mathbf{y}(0) & \in [\mathbf{y}_0], \quad\mathbf{p} \in [\mathbf{p}],\quad t \in[t_0, t_{f}]\enspace
  \end{aligned}
\end{equation}

Equation~\eqref{eq:iivp2} can be seen as the following differential inclusion (using the monotonicity of intervals):

\begin{equation}
\label{eq:difincl}
\begin{aligned}
      \dot{\mathbf{y}}(t;[\mathbf{p}]) & \in \mathbf{f}(t,\mathbf{y}(t;[\mathbf{p}]);[\mathbf{p}])
     \\
      \text{with }\mathbf{y}(0) & \in [\mathbf{y}_0],\quad t \in[t_0, t_{f}]\enspace
  \end{aligned}
\end{equation}

Validated simulation, presented in Section~\ref{sec:validated-numerical-integration-methods}, is able to solve this problem, and provides the functions $R(t)$ and $\widetilde{R}([t])$, for a given $[\mathbf{p}]$. In the following, we will use the notations $R(t;[\mathbf{p}])$ and $\widetilde{R}([t];[\mathbf{p}])$ to keep the dependency of the solution of the IVP with respect to the current parameter.

\subsection{Cost evaluation}
Second, the cost function has to be evaluated with intervals. We recall the cost function of the Problem~\ref{problem}:

 \begin{equation*}
  \min_{\mathbf{p}} J(\mathbf{p}) = \phi(\mathbf{y}(t_f;\mathbf{p}), \mathbf{p}) + \int_{t_0}^{t_f} g(t, \mathbf{y}(t;\mathbf{p}), \mathbf{p}) dt
 \end{equation*}
 
The first part of the cost function, the terminal cost, can be easily bounded with the solution of the IVP:

\begin{equation}
  \phi(\mathbf{y}(t_f;\mathbf{p}), \mathbf{p}) \in \phi(R(t_f;[\mathbf{p}]), [\mathbf{p}])
\end{equation}

The second part, the continuous cost, is evaluated following the discretization using to solve the IVP (the discretization can be increased to obtain a sharper result), and the rectangle rule \cite{Moore66}:

\begin{equation}
 \int_{t_0}^{t_f} g(t, \mathbf{y}(t;\mathbf{p}), \mathbf{p}) dt \in \sum_{i=0}^n (t_{i+1}-t_i) \times g([t_i,t_{i+1}], \widetilde{R}([t_i,t_{i+1}];[\mathbf{p}])], [\mathbf{p}])
\end{equation}

The rectangle rule can be replaced by a more precise scheme such as the Taylor series or the Simpson rule.

\subsection{Branch and Bound algorithm}
\label{alg:bb}
In this section, we present the simplest Interval Branch \& Bound algorithm \cite{JKDW01} that can be used for a reliable minimization of a problem described by:

\begin{equation}
\label{eq:optim_glob}
 \min_p J(\mathbf{p}), \text{ such that } C(\mathbf{p}) \text{ is valid,}
\end{equation}

with $\mathbf{p} \in [\mathbf{p}]$, $J(.)$ an objective function and $C(.)$ a set of constraints (a constraint satisfaction problem for example).  

Many improvements can be considered \cite{hansen_global_2003}, but we prefer the simplest one in order to show the efficiency of our contribution in a clearer way. 

Given an initial box $[\mathbf{p}]$, and a threshold $\epsilon$:

\begin{tabular}{p{2cm} p{11cm}}
 Step 0. & Initialisation\\ 
 & Set the upper bound on the objective function $\Sup{c} = \infty$, the list of optimal solutions $\mathcal{L}=\emptyset$, and the working list of boxes $\mathcal{Q}=\{([\mathbf{p}],\infty)\}$.\\
 Step 1. & Selection\\
 & Pop first box out of $\mathcal{Q}$ into $[\mathbf{p}]$.\\
 Step 2. & Upper bound\\
 & Reduce $[\mathbf{p}]$ with respect to the constraints $C(.)$. If $[\mathbf{p}]$ is empty, then the problem is infeasible in this subregion, go to step 5. Otherwise, compute the objective function in a particular point of $[\mathbf{p}]$ (for example the middle, if feasible), and update the upper bound $\Sup{c} = min(\Sup{J(\Mid{[\mathbf{p}]})} , \Sup{c})$.\\
 Step 3. & Filtering\\
 & Remove from $\mathcal{Q}$ any pair $([\mathbf{p}_i],{c}_i))$ such that $c_i > \Sup{c}$. \\
 Step 4. & Branching\\
 & If $\width{[\mathbf{p}]} < \epsilon$, then put the pair $([\mathbf{p}],\Inf{J([\mathbf{p}])})$ into $\mathcal{L}$. Otherwise, apply a branching rule on $[\mathbf{p}]$ to choose a variable on which to branch and generate two new boxes $[\mathbf{p}_1]$ and $[\mathbf{p}_2]$ (such that $[\mathbf{p}]=[\mathbf{p}_1] \cup [\mathbf{p}_2]$). Put $([\mathbf{p}_1],\Inf{J([\mathbf{p}_1])})$ and $([\mathbf{p}_2],\Inf{J([\mathbf{p}_2])})$ into $\mathcal{Q}$.\\
 Step 5. & Loop\\
 & If $\mathcal{Q}$ is non empty then go to step 1.\\
 Step 6. & Filtering of solutions\\
 & Remove from $\mathcal{L}$ any pair $([\mathbf{p}_i],{c}_i)$ such that $c_i > \Sup{c}$. Compute the union of all remaining solutions in $\mathcal{L}$: $[\mathbf{p}_{sol}] = \bigcup_i [\mathbf{p}_i]$ with a cost $[c_{sol}] = \bigcup_i J([\mathbf{p}_i])$.\\
\end{tabular}

\paragraph*{Step 2: Upper bound\\}
In this paper, we focus on the case of guaranteed global optimization with nonlinear ODEs. The computation of the upper bound of the objective function, done during this step, requires to solve the ODE, that is to say to perform the validated simulation for the current box of parameters. The method presented in Section~\ref{sec:validated-numerical-integration-methods} is used for this purpose.

\paragraph*{Step 4: Branching\\}
In this step, the box $[\mathbf{p}]$ is partitioned into two new boxes. the variable on which to branch is selected following a bisection heuristic. The two main strategies for this purpose are:
\begin{itemize}
 \item Round Robin: each variable is bisected in turn,
 \item Largest First: the splitting is done on the variable having the largest width.
\end{itemize}

The contribution presented in this paper propose a method to improve this step with the help of sensitivity analysis. 

\subsection{Sensitivity analysis}

Sensitivity analysis has been already used in a problem of guaranteed global optimization with ODEs, but in a different manner than the one presented in this paper. Indeed, in \cite{rauh2013interval} authors exploit the sensitivity analysis in order to improve the local search done in step 2. This approach is much more efficient than the trial in the middle as performed in the previous algorithm. For a Branch \& Bound algorithm, it is important to find a ``good'' minimum as soon as possible to prune the branches and in that way limit the number of boxes in the working list. A local search can be added to our method without any loss of generality. 


\section{Contribution: sensitivity function for Smear-based bisection heuristic}
\label{sec:smear}
In this section, the well-known Smear-based heuristic is recalled. Then, our contribution, a similar approach based on sensitivity function, is presented. 

\subsection{Smear-based heuristic}
In the branching-based algorithms, the heuristic of bisection has an important effect on the number of branches, the complexity and finally the rapidity of the solution computing. In the case of root-finding problems (or continuous constraint satisfaction problems), solved with a Branch \& Prune algorithm, it is assumed that the Smear-based heuristic is often the most efficient. This latter is based on the Smear function as defined in \cite{Kearfott:1990:AIP:78928.78931}. 

\begin{definition}[Smear heuristic]
 Let $J=((J_{ij} = \frac{\partial f_i}{\partial x_j}))$ be the Jacobian matrix of the system, and let, for each variable $x_j$, the {\bf smear value} be given by $s_j = {\rm Max}( \vert\underline{J_{ij}[\underline{x_j},\overline{x_j}]}\vert, \vert\overline{J_{ij}[\underline{x_j},\overline{x_j}]}\vert$ $\forall i\in [1,n]$, where $n$ is the total number of functions. 
 The {\bf Smear heuristic} consists to bisect on dimension $k$ if $s_k > s_j, \forall j \in [1,n]\setminus k$.
\end{definition}

{\it According to the article \cite{Kearfott:1990:AIP:78928.78931}, there is a drawback in the use of the smear heuristic. Considering for example the equation $f(a,b)=a^2b^3-1$, where $a$ and $b$ are given in large identical intervals centered at $0$. If we are looking for the zeros of $f$ by the help of a Branch \& Prune algorithm with the Smear heuristic for bisection:

The derivative of $f$ with respect to $a$ is $2ab^3$ and with respect to $b$ $3a^2b^2$: multiplied by the width of the interval we get $2a^2b^3$ and $3a^2b^3$. Hence the smear function for $b$ will be in general larger than for $a$, and $b$ will always be bisected until its width is lower than the desired accuracy.  
However, the smear function is very often the most efficient mode and should be privileged. }

The Smear heuristic has also been used in the case of global optimization \cite{kearfott2013rigorous}. Let consider the problem given in Equation~\eqref{eq:optim_glob}, the Smear leads to bisect on dimension $k$ if: 
\begin{equation}
 \Bigg\vert \frac{\partial J}{\partial p_k}([\mathbf{p}]) \Bigg\vert \width{p_k} = \max_{1\leq i \leq n} \Bigg\vert\frac{\partial J}{\partial p_i}([\mathbf{p}])\Bigg\vert \width{p_i}
\end{equation}

\subsection{Sensitivity function for Smear-based heuristic}
We propose to define a new Smear for the problem of guaranteed global optimization with ODEs by the help of sensitivity function. 
In the problem \ref{problem}, it is important to understand that the objective function $J(\mathbf{p})$ is also function of $\mathbf{y}(t;\mathbf{p})$, the solution of IVP. The objective function is then written $J(\mathbf{p}, \mathbf{y}(t;\mathbf{p}))$, and its derivative with respect to the $i^{th}$ parameter is obtained as follows:

\begin{equation}
 \frac{\partial J(\mathbf{p}, \mathbf{y}(t;\mathbf{p}))}{\partial p_i} = \frac{\partial J}{\partial \mathbf{y}} \frac{\partial \mathbf{y}}{\partial p_i}(\mathbf{p}, \mathbf{y}(t;\mathbf{p}))
\end{equation}

 The quantity $\frac{\partial J}{\partial \mathbf{y}} (\mathbf{p}, \mathbf{y}(t;\mathbf{p}))$ is the same whatever the parameter $p_i$ considered. 
 %
 %
Then, the derivative of objective function with respect to the $i^{th}$ parameter is simplified as follows:

\begin{equation*}
   \frac{\partial J(\mathbf{p}, \mathbf{y}(t;\mathbf{p}))}{\partial p_i} = C \frac{\partial \mathbf{y}}{\partial p_i}(\mathbf{p}, \mathbf{y}(t;\mathbf{p})),
\end{equation*}

$C$ being a function, identical whatever the considered parameter. 
 
\begin{proposition}[Sensitivity function for Smear-based heuristic]
The Smear-based heuristic consists to bisect on $k^{th}$ dimension for which
    \begin{equation}
    \label{eq:smear}
   \Bigg\vert C \frac{\partial \mathbf{y}}{\partial p_k}(\mathbf{p}, \mathbf{y}(t;\mathbf{p})) \Bigg\vert w(p_k) = \max_{1\leq i \leq n} \Bigg\vert C \frac{\partial \mathbf{y}}{\partial p_i}(\mathbf{p}, \mathbf{y}(t;\mathbf{p}))\Bigg\vert w(p_i)
  \end{equation}
  where $\frac{\partial \mathbf{y}}{\partial p_i}(\mathbf{p}, \mathbf{y}(t;\mathbf{p}))$ is the sensitivity function $\mathbf{s}_i(t;\mathbf{p})$ obtained with the method given in Section~\ref{sec:adjoint}. The value $C$ can be obviously simplified. 
  The condition to choose the splitting dimension is then the following:
  \begin{equation}
    \label{eq:smear2}
   \vert \mathbf{s}_k(t;\mathbf{p}) \vert w(p_k) = \max_{1\leq i \leq n} \vert \mathbf{s}_i(t;\mathbf{p})\vert w(p_i)
   \end{equation}
   \hfill$\square$
\end{proposition}

Condition \eqref{eq:smear2} makes appear a norm on the continuous function $\mathbf{s}_i(t; \mathbf{p})$. In the special case of a terminal cost (only function of $t_f$), the problem is simplified. In the other cases, the norm has to be evaluated along time. For that, the discretization used to solve the IVP, and already exploited for the continuous cost, can be considered to compute the norm of the abstracted solution of $\mathbf{s}_i(t; \mathbf{p})$ given by $ \widetilde{R}([t];[\mathbf{p}])]$. 

Using the infinity norm (or supremum):

\begin{equation}
\vert \vert\mathbf{s}_i(t; \mathbf{p}) \vert\vert_{\infty} \leq \max \{ \vert\vert \widetilde{R}([t_j,t_{j+1}];[\mathbf{p}]) \vert\vert_{\infty}\quad ,\quad j=1,\dots,n-1 \}
\end{equation}

This norm is easy and fast to compute, and provides the limit of the $l_p$ norm (while $p \to \infty$) evaluated on the interval $[t_0,t_f]$:

\begin{equation}
 \vert \vert\mathbf{s}_i(t; \mathbf{p}) \vert\vert_{p} = \Big( \int_{t_0}^{t_f} \vert \mathbf{s}_i(t; \mathbf{p}) \vert ^p dt \Big) ^{1/p}
\end{equation}

The infinity (maximum) norm of a box $[\mathbf{a}] \in \IRset^n$ is given by

\begin{equation}
 \vert\vert ~[\mathbf{a}]~ \vert\vert_{\infty} = \max \{\vert~ [a_i]~ \vert \quad , \quad i=1,\dots,n\},
\end{equation}

while the absolute value of an interval $[a_i] \in \IRset$, such that $[a_i]=[\Inf{a_i},\Sup{a_i}]$,  is defined by

\begin{equation}
 \vert ~[a_i]~ \vert = \max \{\vert \Inf{a_i} \vert, \vert \Sup{a_i} \vert\}.
\end{equation}
 
\section{Case studies}
\label{sec:case}
This section examines three example problems. The first one comes from the literature, the well-known singular control problem, originally proposed in \cite{luus1990optimal} and discussed many times in global optimization with non linear ordinary differential equations papers \cite{Singer2006}. This problem considers a one dimensional parameter, thus our contribution cannot be applied, but it allows us to show the efficiency of the global algorithm. 
In order to judge the efficiency of a heuristic, it is common to find a case where the heuristic leads to an important gain, and one where the heuristic does not work (due to the structure of the example), but where the deterioration remains low. Following this method, the second example is proposed to push the limits of our approach (three parameters, a nonlinear objective function and a nonlinear differential equation), and shows its efficiency. The last one, an optimal control problem with a two dimensional control and one end-point constraint, formerly proposed in \cite{goh1988control}, is allows us to produce a counter-example to prove the limits of our heuristic. 


\subsection{Singular control}
\label{ex:sing}
The original problem is given by
\begin{equation*}
 \min_{u} \int_{t_0}^{t_f} y_1^2 +y_2^2+0.0005(y_2+16t-8-0.1y_3u^2)^2 dt
\end{equation*}
subject to
\begin{equation*}
  \left\{
    \begin{aligned}
     \dot{y_1} &= y_2\\
     \dot{y_2} &= -y_3u+16t-8\\
     \dot{y_3} &= u
  \end{aligned}
\right.
\end{equation*}

with $\mathbf{y}(0) = (0,-1,-\sqrt{5})$, $t \in [0,1]$ and $u \in [-4,10]$. The problem is solved with the Algorithm~\ref{alg:bb}, with a tolerance of $1\text{e-}4$. The obtained result is  $u_{sol} \in [3.9003, 4.2165]$ for an upper bound for the cost $c_{sol} < 0.5044$. The algorithm branches 10376 times. The results given in \cite{Singer2006} are $u_{opt}=4.07$ and $c_{opt}=0.497$. The guaranteed results obtained with our unsophisticated algorithm are then compatible with these latters.

\subsection{Polynomial example}
\label{ssec:newprob}
An example is built in order to experience the proposed contribution. The problem is the following:

\begin{equation*}
 \min_{p_1,p_2,p_3} (y_1(t_f) + y_2(t_f))^2
\end{equation*}
subject to
\begin{equation*}
  \left\{
    \begin{aligned}
     \dot{y_1} &= p_1 y_1^2 + p_2 y_2 - 2 p_3^2 \\
     \dot{y_2} &= -3 p_1 y_1 - p_1 p_2 y_2 + y_1 y_2 p_3 +1.0 \\
  \end{aligned}
\right.
\end{equation*}

with $\mathbf{y}(0) = (0,1)$, $t \in [0,1]$ and $p \in [0.95,1]^3$. 

The sensitivity function of the differential equation with respect to $p_1$ is given by the solution of the following ODE:
\begin{equation*}
  \left\{
    \begin{aligned}
    \dot{y_1} &= p_1 y_1^2 + p_2 y_2 - 2 p_3^2 \\
     \dot{y_2} &= -3 p_1 y_1 - p_1 p_2 y_2 + y_1 y_2 p_3 +1.0 \\
    \dot{s_1} &=  2 p_{1} s_{1} y_{1} + p_{2} s_{2} + y_{1}^{2}\\
    \dot{s_2} &= - p_{2} y_{2} + s_{1} \left(- 3 p_{1} + p_{3} y_{2}\right) + s_{2} \left(- p_{1} p_{2} + p_{3} y_{1}\right) - 3 y_{1}~,
    \end{aligned}
  \right.
\end{equation*}

the sensitivity function of the differential equation with respect to $p_2$ is given by the solution of the following ODE:
\begin{equation*}
  \left\{
    \begin{aligned}
    \dot{y_1} &= p_1 y_1^2 + p_2 y_2 - 2 p_3^2 \\
     \dot{y_2} &= -3 p_1 y_1 - p_1 p_2 y_2 + y_1 y_2 p_3 +1.0 \\
    \dot{s_1} &=   2 p_{1} s_{1} y_{1} + p_{2} s_{2} + y_{2}\\
    \dot{s_2} &= - p_{1} y_{2} + s_{1} \left(- 3 p_{1} + p_{3} y_{2}\right) + s_{2} \left(- p_{1} p_{2} + p_{3} y_{1}\right)~,
    \end{aligned}
   \right.
\end{equation*}

and the sensitivity function of the differential equation with respect to $p_3$ is given by the solution of the following ODE:
\begin{equation*}
  \left\{
    \begin{aligned}
    \dot{y_1} &= p_1 y_1^2 + p_2 y_2 - 2 p_3^2 \\
     \dot{y_2} &= -3 p_1 y_1 - p_1 p_2 y_2 + y_1 y_2 p_3 +1.0 \\
    \dot{s_1} &=  2 p_{1} s_{1} y_{1} + p_{2} s_{2} - 4 p_{3}\\
    \dot{s_2} &= s_{1} \left(- 3 p_{1} + p_{3} y_{2}\right) + s_{2} \left(- p_{1} p_{2} + p_{3} y_{1}\right) + y_{1} y_{2}~.
    \end{aligned}
   \right.
\end{equation*}

Two bisection strategies are used, the largest first and the sensitivity-based Smear. Results are gathered in Table~\ref{tab:newprob}. The theoretical solution being $(0.95,1,1)$, our results are then consistent.


 \begin{table}
 \begin{tabular}{p{1cm} | p{5cm} p{1.2cm} p{1.2cm} p{1.2cm} p{1cm}}
 Prec. & Optimal solution & Cost & Branches (LF) & Branches (S) & Gain\\
 \hline
$10^{-2}$ & $([0.95, 1] ; [0.95, 1] ; [0.969, 1])$ & 0.72805 & 399 & 357 & 11\%\\
$10^{-3}$  & $([0.95, 0.977] ; [0.992, 1] ; [0.996, 1])$& 0.71991 & 4331 & 3154 & 27\%\\
$10^{-4}$  & $([0.95, 0.954] ; [0.998, 1] ; [0.998, 1])$ & 0.71889 & 8999 & 6460 & 28\%\\
$10^{-5}$  & $([0.95, 0.951] ; [0.999, 1] ; [0.999, 1])$& 0.71875 & 16864 & 12154 & 28\%\\
 \end{tabular}
 \caption{Results of our algorithm on Problem~\ref{ssec:newprob} with Largest First (LF) and Smear (S) strategies.}
 \label{tab:newprob}
\end{table}

\subsection{Optimal control with end-point constraint}
\label{ssec:endcont}
The formulation of the problem is:

\begin{equation*}
 \min_{u_1,u_2} y_2(t_f) 
\end{equation*}
subject to
\begin{equation*}
  \left\{
    \begin{aligned}
     \dot{y_1} &= u_1(1-t)+u_2t\\
     \dot{y_2} &= y_1^2 + (u_1(1-t)+u_2t)^2)\\
  \end{aligned}
\right.
\end{equation*}

with $\mathbf{y}(0) = (1,0)$, $t \in [0,1]$ and $u \in [-1,1]^2$. A constraint on the end-point is given by $y_1(t_f) = 1$. 
 This problem has been previously studied with a rigorous approach in \cite{papamichail2002rigorous}. In this latter, the results are $u_1=-0.4545$, $u_2=0.4545$, with a value of the objective function for the global optimum parameters equal to $0.924242$. 

The sensitivity function of the differential equation with respect to $u_1$ is given by the solution of the following ODE:
\begin{equation*}
  \left\{
    \begin{aligned}
     \dot{y_1} &= u_1(1-t)+u_2t\\
     \dot{y_2} &= y_1^2 + (u_1(1-t)+u_2t)^2)\\
     \dot{s_1} &= - t + 1\\
     \dot{s_2} &= 2 s_{1} y_{1} + \left(- 2 t + 2\right) \left(t u_{2} + u_{1} \left(- t + 1\right)\right)\\
  \end{aligned}
\right.
\end{equation*}

While the sensitivity function of the differential equation with respect to $u_2$ is given by the solution of the following ODE:
\begin{equation*}
  \left\{
    \begin{aligned}
     \dot{y_1} &= u_1(1-t)+u_2t\\
     \dot{y_2} &= y_1^2 + (u_1(1-t)+u_2t)^2)\\
     \dot{s_1} &= t\\
     \dot{s_2} &= 2 s_{1} y_{1} + 2 t \left(t u_{2} + u_{1} \left(- t + 1\right)\right)\\
  \end{aligned}
\right.
\end{equation*}

The results obtained are given in Table~\ref{tab:endcont}. With Algorithm~\ref{alg:bb}, the tightest global optimum (obtained with $\epsilon=\num{1e-5}$) is such that 
\[
(u_1,u_2) \in ([-0.462448, -0.446724] ; [0.446716, 0.46244])
\]
 , with an upper bound on the objective function valued to $0.924249$. Our results are then consistent with the literature. 

\begin{table}
  \begin{tabular}{p{1cm} | p{5cm} p{1.2cm} p{1.2cm} p{1.2cm} p{1cm}}
 Prec. & Optimal solution & Cost & Branches (LF) & Branches (S) & Loss\\
 \hline
$10^{-3}$ & $([-0.511, -0.399] ; [0.398, 0.510])$ & 0.924249 & 1243 & 1250 & 0.5\%\\
$10^{-4}$ & $([-0.469, -0.439] ; [0.439, 0.469])$ & 0.924249 & 5362 & 5369 & 0.2\%\\
$10^{-5}$ & $([-0.462, -0.446] ; [0.446, 0.462])$ & 0.924249 & 19283 & 19290 & 0.1\%\\
 \end{tabular}
 \caption{Results of our algorithm on Problem~\ref{ssec:endcont} with Largest First (LF) and Smear (S) strategies.}
 \label{tab:endcont}
\end{table}




\subsection{Discussion}
Through experiments, we demonstrate that our algorithm is able to solve a guaranteed global optimization with ODEs as shown in the first example~\ref{ex:sing}. In the second example~\ref{ssec:newprob}, we exhibit the efficiency of the proposed heuristic with a gain up to 30\% in term of number of branches. A third example~\ref{ssec:endcont}, where our heuristic fails to improve the solving procedure, leads to a very small loss, below 0.5\%.

\section{Conclusion}
\label{sec:concl}
Differential equations are widely used to model physical, chemical and biological systems. However, there has been little development in term of guaranteed solvers - solvers which provide validated solutions with respect to method approximations, uncertainty on parameters and rounding error - for the most important problems such as parameter identification, optimization, constraint satisfaction problems, etc. The branching based algorithms (with interval arithmetic) are often the most efficient approaches for this class of problems. However, the number of branches has to be reduced as much as possible. In this paper, a new heuristic for bisection is proposed for the special case of optimization involving ODEs. For this purpose, the sensitivity function is computed to provide a version of the Smear strategy. 

The obtained heuristic is used in a basic branch-and-bound algorithm to solve some examples. If in some cases the heuristic does not reduce the number of branching, sometimes it can reduce by thirty percents the number of branches in the search tree. 

The bisection strategy proposed can be used in a more sophisticated branch-and-bound algorithm. As future work, it can be interesting to reuse the sensitivity function computed in a local search approach to accelerate the upper bound decreasing. In term of implementation, some efforts are needed to accelerate the computation of the sensitivity function. 

To conclude, a novel heuristic has been proposed to reduce the number of branches required to find the optimal in a branch-and-bound algorithm. The results are encouraging.

\section*{Softwares}
All the results in this paper have been obtained with DynIbex, a tool for validated simulation, and with Sympy for the symbolic computation of sensitivity functions. 

\section*{Nomenclature/Notation}
$\dot{y}$ denotes the time derivative of $y$, that is, $\frac{d y}{d
  t}$. $a$ denotes a real value, while $\textbf{a}$ represents a vector of real values. $[a]$ represents an interval value and $[\mathbf{a}]$ represents a vector of interval values (a box). The midpoint of an interval $[x]$ is denoted by $\Mid([x])$. The variables $y$ are used for the state variables of the system and $t$ represents time. Sets will be represented by calligraphic letter such as $\mathcal{X}$ or $\mathcal{Y}$. An interval with floating point bounds is written in the short form, for example $0.123456[7,8]$ to represent the interval $[0.1234567,0.1234568]$.

\paragraph{Funding:}
This research benefited from the support of the ``Chair Complex Systems Engineering - Ecole Polytechnique, THALES, DGA, FX, DASSAULT AVIATION, DCNS Research, ENSTA ParisTech, T\'el\'ecom ParisTech, Fondation ParisTech and FDO ENSTA''. This research is also partially funded by DGA MRIS.

\paragraph{Conflict of Interest:} 
The authors declare that they have no conflict of interest.


\bibliographystyle{plain}
\bibliography{bib}

\end{document}